\documentclass[12pt]{article}

\usepackage{amsmath,amssymb,amsthm,amscd,amsfonts,float,graphicx}
\usepackage[applemac]{inputenc}
\usepackage[all]{xy}
\theoremstyle{plain}
\newtheorem{theor}{Theorem}[section]

\newtheorem{Prop}[theor]{Proposition}

\theoremstyle{definition}
\newtheorem{defin}[theor]{Definition}

\newcommand{\NN}{(N_t)_{0\leq t \leq 1} }
\newcommand{\ust}{b}
\newcommand{\half}{ {\scriptstyle{\frac{1}{2}} } }
\newcommand{\N}{\mathbb N}

\newcommand{\E}{\mathrm E}

\newcommand{\T}{{\cal T}}

\newcommand{\F}{{\cal F}}
\newcommand{\sF}{\F}
\newcommand{\sT}{\T}

\usepackage{url}
\date{}

\usepackage[authoryear]{natbib}
\begin{document}
\title{Answer to an open question concerning the  $1/e$-strategy for best choice under no information.}
\author{F. Thomas Bruss\footnote{F.\,Thomas Bruss, Universit\'e Libre de Bruxelles, 
D\'epartement de Math\'ematique, CP 210, B-1050 Brussels, Belgium (tbruss@ulb.ac.be)} and L.C.G. Rogers\footnote{L.C.G. Rogers, Statistical Laboratory, University of Cambridge, Wilberforce Road, Cambridge CB3 0WB,
United Kingdom (lcgr1@cam.ac.uk)}\\Universit\'e Libre de Bruxelles and University of Cambridge}
\maketitle
\centerline{\bf In Memory of} \centerline {\bf Professor Larry Shepp}

\bigskip\noindent
\begin{abstract}
 This paper answers a long-standing open question concerning the   $1/e$-strategy for the problem of best choice. $N$ candidates for a job arrive at times independently uniformly distributed in $[0,1]$. The interviewer knows how each candidate ranks relative to all others seen so far, and must immediately appoint or reject each candidate as they arrive. The aim is to choose the best overall. The $1/e$ strategy is to follow the rule: `Do nothing until time $1/e$, then appoint the first candidate thereafter who is best so far (if any).'
 
 The question, first discussed with Larry Shepp in 1983, was to know whether the $1/e$-strategy is optimal if one has `no information about the total number of options'. Quite what this might mean is open to various interpretations, but we shall take the proportional-increment process formulation of \cite{BY}. Such processes are shown to have a very rigid structure, being time-changed {\em pure birth processes}, and this allows some precise distributional calculations, from which we deduce that the $1/e$-strategy is in fact not optimal.
 
\end{abstract}


\medskip\noindent{\bf MSC 2010 Subject Code}: 60G40

\section{Dedication and background}
At the evening of Professor  Larry Shepp's talk ``Reflecting Brownian Motion" at Cornell University on July 11, 1983 (13th Conference on Stochastic Processes and Applications), Professor Shepp and Thomas Bruss ran into each other in front of the Ezra Cornell statue. Thomas was   honoured to meet Prof. Shepp in person, but Larry replied ``What are you working on?" And so Larry was the very first person with whom Thomas could discuss the {\it $1/e$-law of best choice} resulting from the {\it Unified Approach}  \cite{B84} which had been accepted for publication shortly before. 
Thomas was pleased to see the true interest Prof. Shepp showed for the $1/e$-law. As many of us  have seen before, when Larry was interested in a problem, elementary or not, then he was really deeply interested.\smallskip

This article deals with an open question concerning  the so-called $1/e$-strategy for the problem of best choice, which is to wait and do nothing until time $1/e$ and then to accept  the first candidate (if any) who is best so far.  The question which attracted our particular interest was whether this strategy is optimal if one has no initial information about the number $N$ of candidates. Bruss drew also attention to this open question in his own talk ``The $e^{-1}$-law in best choice problems" at Cornell on July 14, 1983, and re-discussed it with Larry at several later occasions. 
A written record of somewhat related questions  appeared on page 885 of \cite{B84} where he stated the conjecture that the $1/e$-strategy is optimal in certain two-person games for a decision maker who faces an adversary trying to minimize the win probability. However, the two-person game situation is quite different from the open question discussed with  Larry and will not be considered in this paper.

\smallskip
 As far as we are aware, the last time the question discussed with Larry was addressed was in \cite{BY}, and this may actually be the only written reference to the real open question. \cite{BY} studied another no-information stopping problem, the so-called last-arrival problem (l.a.p.). To prepare the paper's main result, they examined the hypothesis of no-information  in a detailed way, and their conclusions will be used in the present paper. \cite{BY} also used these to give, as a side-result, an alternative proof of the $1/e$-law. However, as they pointed out, their approach did not contribute new insights for the open question.

\medskip
The present article proves that the $1/e$-strategy is {\it not} optimal under the interpretation of `no information' used by \cite{BY}. It thus closes a 37-years gap.

\section{The Unified Approach}
The so-called {\it 1/e-law of best choice~} is a result obtained in the {\it Unified Approach}-model of \cite{B84}. The model is as follows:
 \begin{quote}{\bf Unified Approach}: Suppose $N>0$ points are IID  $U[0,1].$
Points are marked with qualities which are supposed to be uniquely rankable from $1$ (best) to $N$ (worst), and all rank arrival orders are supposed to be equally likely. The goal is  to maximize the probability of stopping online, and without recall on a preceding observation, on rank 1. \end{quote}
\noindent This model  was suggested for the best choice problem (secretary problem) for an unknown number $N$ of candidates. (More general payoff-functions for the same model were studied in \cite{BrussSamuels}.) 

\medskip


Now, if we contemplate the probability of picking the best candidate, we immediately face the question `What is $N$?'  If $N$ is fixed and known, this is just the classical secretary problem. But if we take a Bayesian point of view and suppose a prior distribution for $N$, with arrivals coming at times $t \in \N$,
 \cite{abdel} showed that the problem may not only lead to so-called {\it stopping islands} (\cite{presman}), but, much worse, that for any $\epsilon>0$ there exists a sufficiently unfavorable distribution $\{P(N=n)\}_{n=1,2, \cdots}$ to reduce the value of the optimal to less than   $\epsilon.$ In other words, if $N$ is allowed an arbitrary prior,  optimality may mean almost nothing.
These discouraging facts  prompted efforts to find more tractable models, such as the model of \cite{stewart}, and the one of \cite{cowan} and its generalisation studied in \cite{bruss88}.

\smallskip The philosophy behind the unified approach of \cite{B84} was different. The approach was to suppose that arrival times are in $[0,1]$ and to study so-called $x$-strategies, where you do nothing until time $x$, and thereafter pick the first record. One of the main results of that paper was that the $1/e$-strategy gives a success probability of at least $1/e$ whatever the prior distribution of $N$, and that no other $x$-strategy does this well. This robustness suggests that the $1/e$-strategy is somehow special, and the open question became natural.

It is relevant to mention here that a similar phenomenon of robustness shows up in different forms. One is what \cite{bruss1990conditions} called `quasi-stationarity', meaning essentially that the optimal strategy
may (even for rather general payoffs) hardly depend on the number of candidates observed so far.
More remarkably, for so-called Pascal processes, optimal strategies do not depend at all on the number of preceding observations (For their characterization see \cite{bruss1991pascal}). 

\subsection{The open question}\label{conj}
First we need to be clear about what exactly we mean by optimality of a strategy under no information on $N.$  We see a counting process $\NN$, $N_0=0$, and we define $\sF_t = \sigma( N_u,u \leq t)$. The law of $\NN$ is $P_\theta$ for some $\theta \in \Theta$, where $\{P_\theta:  \theta \in \Theta\}$ is the collection of possible laws considered. 

The notion that `we have no prior information at all on $N$' means that we are only going to consider strategies which are $(\sF_t)$-stopping times.  That is, the strategies allowed can only know the arrival times (and ranks) of the individuals, not the value of $\theta \in \Theta$. This is the viewpoint of classical statistics.

\medskip\noindent
To understand the sense of optimality, define the process $\rho$ by
\begin{eqnarray*}
\rho_t &=& \hbox{\rm overall rank of object arriving at $t$ if $\Delta N_t =1$}
\\
&=& 0 \quad \hbox{\rm otherwise.}
\end{eqnarray*}
Let $\sT$ denote the set of all $(\sF_t)$-stopping times. Then the value of using $\tau \in \sT$ is 
\begin{equation}
R(\theta,\tau) = P_\theta[ \rho_\tau = 1 ].
\label{Rdef}
\end{equation}
We denote by $\tau^*$ the stopping time corresponding to the $1/e$ strategy, which is  simply $\tau^* = \inf \{ t \geq 1/e:  \rho_t = 1 \},$ where it is understood that $\tau^*=1$ if no such $t$ exists, and that, in this case, we lose by definition. In these terms, the open question is stated precisely as follows:  

True or false
\begin{equation}
\forall \theta \in \Theta, \forall \tau \in \sT, \qquad
R(\theta, \tau^*) \geq R(\theta, \tau)?
\label{1/e_conj}
\end{equation}
Of course, the set $\Theta$ of possible laws of $\NN$ plays an important r\^ole in the conjecture. For example, if $\Theta$ contained just one law, under which $\NN$ was the counting process of ten $U[0,1]$ arrival times, then clearly the $1/e$ strategy would not be optimal in the sense of \eqref{1/e_conj}. We shall shortly explain exactly what set of laws is considered here.

\subsection{A related problem.}\label{ss22}

\smallskip
We return to the related {\em last-arrival problem}
under no information (l.a.p.) studied in \cite{BY}.
In this model an unknown number $N$ of points are  IID $U[0,1]$ random variables, and an observer, inspecting the interval $[0,1]$ sequentially from left to right, wants to maximise the probability of stopping online on the very last point. No information about $N$ whatsoever is given.
Only one stop is allowed, and this again without recall on  preceding observations.

\medskip 

Central to the approach of \cite{BY} is the choice of the family $\Theta$ of laws of the counting process $\NN$. These authors present arguments (based on the properties of IID $U[0,1]$ arrival times) to justify their focus on the family of what they call {\em proportional-increments (p.i.)} counting processes. We shall not repeat all the reasoning which leads to this choice of counting processes, but we show its basic motivation and explain why we take its implications as our starting point. 

\smallskip
\cite{BY}
defined a p.i.-process as follows:
A stochastic process $(N_t)$ defined on a filtered probability space$(\Omega, {\cal F}, ({\cal F}_t), P)$ with natural filtration ${\cal F}_t=\sigma\{N_u: u\le 1\}$ is a p.i.-counting process on $]0,\infty[$, if $$\forall t ~{\rm with}~ N_t>0, ~\forall s \ge 0,$$
$$\E(N_{t+s}-N_t \,\Big |\, {\cal F}_t) = \frac{s}{t} N_t~a.s.$$ The meaning of {\it proportional} is well understood in this definition. Moreover, three out of the four conclusions 1.- 4. in \cite{BY}, implying this definition, are proved to be compelling for combining the IID $U[0,1]$ - hypothesis for arrival times with the hypothesis that no prior information on $N$ can be used. Only Conclusion 3 (on page 3244) makes a concession. Here these authors use an (unprovable) tractability argument to justify why an unknown random variable, of which the expectation must be zero, is put equal to zero. 

Why a concession? It is important to note that, if one has no information on $N,$ then the time of the first arrival $T_1$ is a particularly delicate point. It is the smallest order statistic of all arrival $N$ times. However, it is exactly this one which escapes any distributional prescription because the no-information setting does not allow us to assume a prior distribution $\{P(N=n)\}_{n=1,2, \cdots}.$ 
Hence, if one wants to confine one's interest to a well-posed problem, as \cite{BY} did, one has to make a concession somewhere if one wants to properly define a relevant decision process in the no-information case. 
The mentioned concession seemed the least restrictive and almost compelling, but, more importantly, \cite{BY} found a solid a-posteriori justification for their tractability argument. The solution of the l.a.p. they obtained for p.i.-processes
satisfied the  criteria of \cite{hadamard}
for the solution of a well-posed problem.  \cite{BY} found these criteria convincing.

Now note that the only difference between the l.a.p. and our open problem (how to find rank 1) is that we want to stop on the last record of the arrival process, and not on the last point.  By the IID-hypothesis for arrival times of absolute ranks, R\'enyi's Theorem of relative ranks (\cite{renyi}) implies that the $k$th point is a record with probability $1/k$ independently of preceding arrivals. Thus the basic arrival process
$(N_t)$ is not affected and can be chosen exactly the same!
 This is why, confining our interest to  well-defined problems only, we suppose that $(N_t)$ is a p.i.-process in the sense of \cite{BY}, from which we take the following definition.
\medskip

\begin{defin}
{\em A p.i.- counting process is a counting process whose compensator is $\lambda_t \equiv N_t/t$, so that ($t\in (0,1]$)
\begin{equation}
M_t \equiv N_{t \vee T_1} - N_{T_1} - \int_{T_1}^{t \vee T_1} \frac{N_s}{s}\; ds
\qquad \hbox{\rm is a martingale in its own filtration,}
\label{PIdef}
\end{equation}
where $T_1 \equiv \inf\{ t: N_t=1\}$ is the first jump time of the counting process.}
\end{defin}

\medskip
 The class $\Theta$ of counting processes will be the class of all p.i.-processes, and the meaning of all the notation appearing in the statement \eqref{1/e_conj} has now been defined.

\section{Analysis of the open question.}\label{S3a}
Our analysis starts with the following little result, whose proof is immediate from the statement.

\begin{Prop}\label{prop1}
Suppose that $(N)$ is a p.i.-counting process.
If we define $\tilde{N}(u) = N(e^u)$ for $u \in (-\infty,0]$, and $t_1
= \log T_1$,  then
\begin{eqnarray*}
M(e^u) &=& N(e^u\vee T_1)-N(T_1) - \int_{T_1}^{e^u \vee T_1} \frac{N_s}{s}\; ds
\\
&=& \tilde{N}(u \vee t_1) - \tilde{N}(t_1) - \int_{t_1}^{u \vee t_1} \tilde{N}(s)
\; ds
\end{eqnarray*}
is a martingale in its own filtration, so  $(\tilde{N})$ is a pure birth process, started with one individual at time $t_1$.
\end{Prop}

\medskip

So the requirement that $(N_t)$ be a p.i.- counting process is not in fact very general - apart from the choice of the time $(\tilde{N})$ starts, the behaviour is uniquely determined!  

\medskip
\noindent
{\sc Remarks.}  If we model a Poisson process with intensity $\lambda$ in a Bayesian fashion, suppose a prior density $f(\lambda) = \varepsilon \exp(- \varepsilon \lambda)$ for $\lambda$, then the posterior mean for $\lambda$ given $\sF_t$ is $N_t/t(1+\varepsilon)$, so a PI counting process is in some sense a limit of a Poisson process where we put an uninformative prior over $\lambda$.

\bigbreak
If we run a pure birth process from $\tilde{N}_u=1$ ($u <0$) to time $0$, the PGF 
of $\tilde{N}_0$ is easily shown to be 
\begin{equation}
E[z^{\tilde{N}_0} \vert \tilde{N}_u = 1] = \frac{ze^{u}}{1-z(1-e^{u})} 
\qquad (z \in [0,1]),
\label{pbp1}
\end{equation}
so that $\tilde{N}_0$ is 1+geometric($e^{u}$). 
Obviously, from \eqref{pbp1} we deduce
\begin{equation}
E[z^{\tilde{N}_0} \vert \tilde{N}_u = k] = \biggl\lbrace
\frac{ze^{u}}{1-z(1-e^{u})} \biggr\rbrace^k
\qquad (z \in [0,1]),
\label{pbp2}
\end{equation}
 Thus if we see a record in the process  $(\tilde{N})$ at time $u <0$, at the arrival of the $n^{ \hbox{\rm th} }$ 
observation, the probability that this is the best overall will be
\begin{equation}
\tilde{\pi}_n(u) \equiv E \biggl[ \;\frac{n}{\tilde{N}_0}\;
\biggl\vert\; \tilde{N}_u = n  \biggr]
= n \int_0^1 \frac{dz}{z} \biggl\lbrace
\frac{ze^{u}}{1-z(1-e^{u})} \biggr\rbrace^n .
\label{pit}
\end{equation}
In terms of the original process $(N)$, if we see a record at the arrival of the 
$n^{ \hbox{\rm th} }$  observation at time $t \in (0,1)$, then the probability that
this is the best overall is
\begin{equation}
\pi_n(t) \equiv E \biggl[ \;\frac{n}{N_1}\;
\biggl\vert\; N_t = n  \biggr]
= n \int_0^1 \frac{dz}{z} \biggl\lbrace
\frac{z t}{1-z(1-t)} \biggr\rbrace^n .
\label{pin}
\end{equation}
Clearly this has to be increasing in $t$, but from numerics it appears also to be decreasing in $n.$ 
We can prove that this has to be the case, as follows. If we fix $t \in (0,1)$ then conditional on $N_t=n$ we have that
\begin{equation}
\frac{N_1}{n} =  \xi_n \equiv\frac{n + W_1+ \ldots + W_n}{n},
\end{equation}
where the $W_j$ are IID geometrics. Now $(\xi_n)$ is a reversed martingale in the exchangeable filtration, so $(1/\xi_n)$ is a reversed submartingale in the exchangeable filtration, so its expectation decreases with $n$.

\section{The value of a fixed threshold  rule.}\label{one_over_e}
Suppose we use a fixed threshold  rule, that is, we do nothing until $u\geq \ust$ and then we take the first record thereafter. The $1/e$ rule corresponds to the special case $\ust=-1$. What is the value of this?

\medskip
If $\tilde{N}_{\ust} = n$, then the distribution of the number $Y$ of further observations is known, and is a negative binomial distribution:
\begin{equation}
P[ Y = y ]  = q^y p^n \binom{n+y-1}{y} \qquad(y \geq 0),
\label{NBdist}
\end{equation}
where $p  = \exp(\ust)$. Given that $Y=y$, the probability that the best comes after the first $n$ observations is $y/(n+y)$, and the probability that the first record after $u=b$ is actually the best is
\begin{equation}
P[ \hbox{\rm first record after $n$ is best}| \hbox{\rm best comes after first $n$,
$Y=y$}]
= \frac{1}{y} \sum_{j=1}^y \frac{n}{n+j-1} \;.
\end{equation}
Thus we have an expression for the probability that we pick the best using this rule:
\begin{equation}
P[ \hbox{\rm win}] = \sum_{y \geq 1} P[Y=y]\;  \frac{n}{n+y}\; \sum_{j=1}^y\;
\frac{1}{n+j-1}.
\label{winprob}
\end{equation}
\subsection{The special case $n=1$.}\label{ss1}
Let us firstly observe that for $t \in (0,1)$
\begin{equation}
f_j(t) \equiv \sum_{k \geq j} \frac{t^k}{k} = \int_0^t \; \frac{s^{j-1}}{1-s} \; ds,
\label{fjdef}
\end{equation}
from which we see that $f_1(t) = -\log(1-t)$.
In the special case $n=1$, we have 
\begin{eqnarray}
P[ \hbox{\rm win}] &=& \sum_{k \geq 1} q^k p\;  \frac{1}{1+k}\; \sum_{j=1}^k\;
\frac{1}{j}\nonumber
\\
&=& pq^{-1} \sum_{j \geq 1} \; \frac{1}{j} \;  f_{j+1}(q)\nonumber
\\
&=& pq^{-1} \int_0^q \sum_{j\geq 1} \frac{s^j}{j}  \; \frac{ds}{1-s} \nonumber
\\
&=& pq^{-1} \int_0^q \biggl( \; \int_0^s \frac{dv}{1-v}
 \; \biggr)  \; \frac{ds}{1-s} \nonumber
 \\
 &=& \half\,  pq^{-1} \biggl( \; \int_0^q \frac{dv}{1-v}
 \; \biggr)^2 \nonumber
 \\
 &=& \half\,  pq^{-1} \bigl( \; \log(1-q)
 \; \bigr)^2 .\label{V1}
\end{eqnarray}
Similarly, from \eqref{pit} we have ($p \equiv 1-q \equiv e^u$)
\begin{eqnarray}
\tilde{\pi}_1(u) &=& \sum_{k \geq 0}\frac{q^k p}{1+k} \nonumber
\\
&=& \frac{p}{q} \; f_1(q)  \nonumber
\\
&=& - \frac{p}{q}\;  \log(1-q).
\label{pi1}
\end{eqnarray}

\section{The Hamilton-Jacobi-Bellman equations.}\label{S3}
If $V_n(u)$ denotes the value of being at time $u\leq 0$ with $n$ events already observed, none of them at time $u$, then the HJB equations of optimal control for the $V_n$ are
\begin{eqnarray}
0 &=& \dot{V_n}(u) + n \biggl\lbrace \frac{n}{n+1} \, V_{n+1}(u) +
\frac{1}{n+1} \max \{ V_{n+1}(u), \tilde{\pi}_{n+1}(u) \} -V_n(u) \; \biggr\rbrace
\nonumber
\\
&=&\dot{V_n} + n( V_{n+1} - V_n ) + \frac{n}{n+1} (\tilde{\pi}_{n+1} - 
V_{n+1})^+,
\label{HJB}
\end{eqnarray}
together with the boundary conditions $V_n(0)=0$. 
The solution is then the value function $V_n(u).$ 
\smallskip

\begin{figure}[H] 
  \centering
   \includegraphics[scale=0.3]{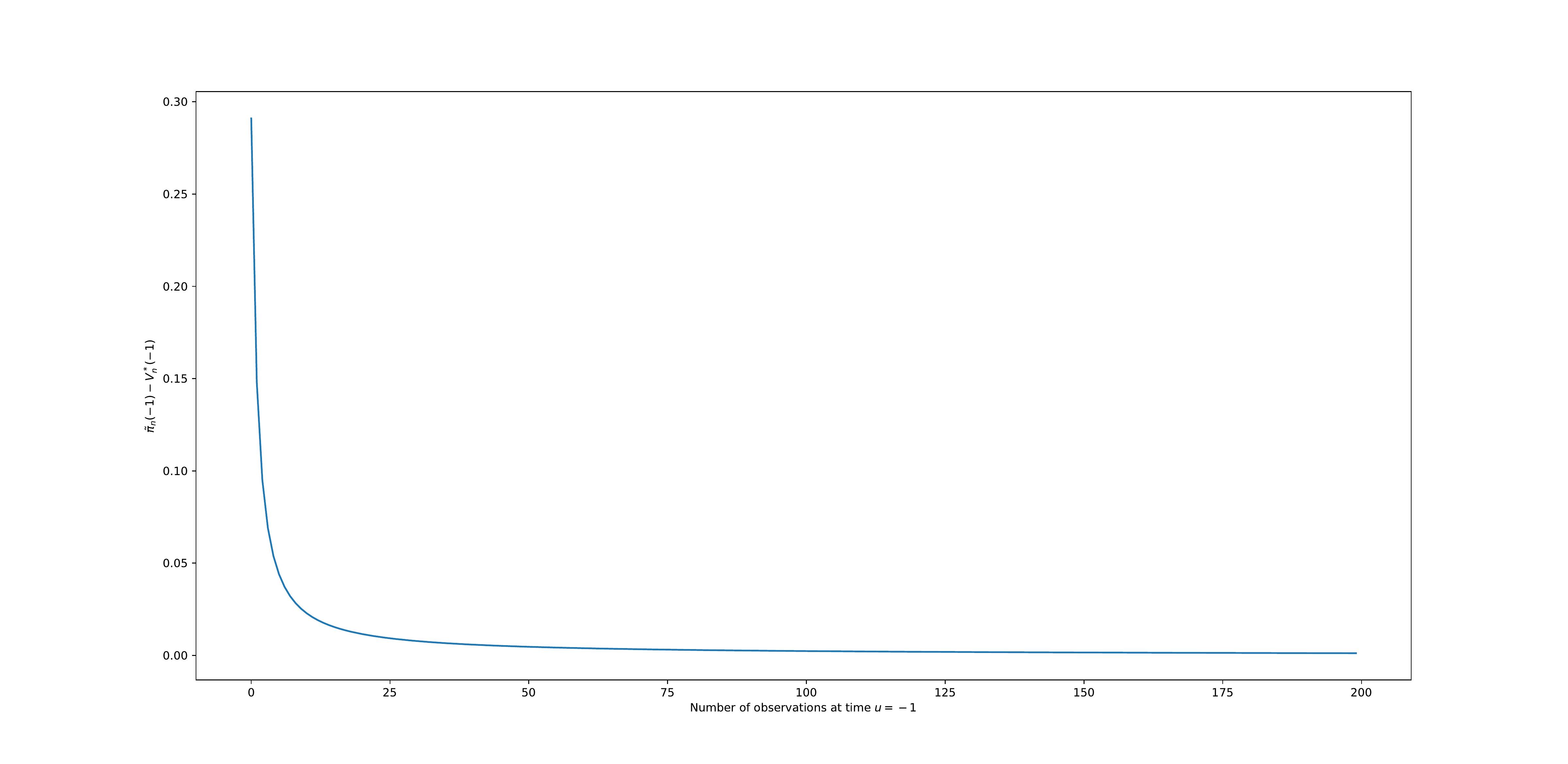} 
   \caption{$\tilde{\pi}_n(-1) - V^*_n(-1)$.}
   \label{ceg}
\end{figure}

\medskip
Now, if the answer to the open question is affirmative, then for all $n$,  $V_n=V_n^*$, where $V^*$ is the value function for using the $1/e$ strategy, which we actually know reasonably explicitly; it is given by the right-hand side of \eqref{winprob}, where the dependence on the time $u<0$ comes via the parameter $p = e^u$ of the negative binomial distribution \eqref{NBdist}.  Now by inspection of the HJB equations, we see that the optimal rule will be that we stop when we see a new record at time $u<0$ when there are $N_u = n$ observations in total if and only if
\begin{equation}
   \tilde{\pi}_n(u) > V_n(u).
   \label{optstop}
\end{equation}
{\em If the $1/e$-strategy is optimal, this would say that}
\begin{equation}
\tilde{\pi}_n(u) > V^*_n(u) \qquad \hbox{\rm if and only if}\quad u > -1.
\end{equation}
We can investigate this numerically by calculating $\tilde{\pi}_n(-1)$ and $V^*_n(-1)$ and comparing them; the results are plotted in Figure \ref{ceg}. 
\noindent
We see that $\tilde{\pi}_n(-1) > V^*_n(-1)$ for all $n$, and that the gap narrows as $n$ increases. Hence the answer to the open question is No. The $1/e$-strategy is not optimal. 
\subsection{Analytic proof.}\label{ss2}
It is nice to see without resorting to numerics that the answer must be No, by considering the special case $n=1$. From \eqref{V1} and \eqref{pi1} we see that for the $1/e$ rule where $p = e^{-1}$, $u=-1$ 
\begin{equation}
\tilde{\pi}_1(u) - V_1(u) = \frac{p}{2q}\bigl[\; -2\log(1-q) - 
\bigl( \; \log(1-q)
 \; \bigr)^2
\;\bigr] = p/2q
\label{gap}
\end{equation}
which is clearly positive.

\subsection{Should this have been obvious?}
A closer look at \eqref{gap} may attract our attention; we may be surprised how big the difference between
$\tilde \pi_1(u)$ and the performance of the $1/e$-strategy can actually become, namely the first one is twice the second one for $u=-1$ (i.e. for $t=1/e$ in $[0,1]$-time.) 
Is it then not surprising that the non-optimality of the $1/e$-strategy did not follow already from simpler comparisons?

No, not easily. First note that we compare here a conditional win probability, given that we have a (granted) record at time $T_1\le 1/e$,  with the unconditional win probability of the 1/e-strategy. Fortunately, this was all that was needed to show that the $1/e$-strategy is not optimal, namely showing that there are situations where it is definitely strictly sub-optimal
to pass over the first arrival, even when arriving at some time $1/e-\epsilon$ with $\epsilon >0$. 

This however does, a priori,  not say much about the absolute win probability of the $1/e$-strategy! To see this, suppose that for some small $\epsilon>$ we have $T_1\in[1/e-\epsilon,1/e+\epsilon]$, and that $N$ is not large. The latter is quite probable if $T_1$ is close to $1/e.$ Then $T_1$ will be almost equiprobable in the left or right half of $[1/e-\epsilon,1/e+\epsilon]$ for $\epsilon$ sufficiently small. If it is in the right half, however, the $1/e$-strategy
will accept it all the same, but now with a strictly greater win probability since $\tilde \pi_1(t)$ is strictly increasing in $t$.  

A second reason why
non-optimality is not that evident lies in the interplay of time and the number of arrivals (see (6) and (8)).
If we have (by simpler estimates) no sufficient incentive to accept at time $t<1/e$ a record, being the $n$th arrival, say, we have even less incentive if it was the $(n+1)$th arrival. With each additional arrival before time $1/e$ increases the expected number of arrivals thereafter, and so does the expected number of those arriving after time $1/e.$ This increases the incentive to wait, but then we also have to wait for a record! This may bring us quickly behind
time $1/e.$ If we get there a record, then, as said in the first scenario, so much the better.

\smallskip
These two scenarios exemplify how important it is to have  precise answers, and not only estimates, even if one has reasonably good ones.
It is the preceding approach which offered these precise answers, and it made things definite: 
In the no-information case, the $1/e$-strategy is {\it not} optimal.

\bigskip \bigskip
\textwidth 14 cm
{\bf Acknowledgement} 

\medskip
The authors would like to thank very warmly Professor Philip Ernst for his  precious indirect and direct contributions to this article. When Philip organised in 2018 the most memorable conference at Rice University in honour of Larry Shepp, he motivated many of us to look at harder problems. And so this old open problem turned up again and attracted attention. Philip's interest in this  problem and his numerous interesting comments on, and discussions of, a former version of this paper, were very encouraging.

\bigskip\bigskip

\bibliography{TBCR}

\end{document}